\documentclass[12pt]{article}
\usepackage [english]{babel}
\usepackage {amssymb}
\usepackage {amsmath}
\begin{document}

\vspace*{10mm}

\large

\begin{center}
{\Large{
Inverse periodic shadowing properties.
\\

\bigskip

A.V. Osipov\footnote{The research of the author is supported by the Chebyshev Laboratory
(Department of Mathematics and Mechanics, Saint-Petersburg State
University) under the grant 11.G34.31.0026 of the Government of the
Russian Federation 
and by Leonhard Euler program financed by German Academic Exchange Service (DAAD). 
}}}
\\

\bigskip

14th Line, 29b, Saint-Petersburg, 199178 Russia\\
Chebyshev Laboratory, Department of Mathematics and
Mechanics, Saint-Petersburg State University\\
e-mail: osipovav@list.ru

\end{center}

\bigskip

\begin{center}
\textbf{Abstract.} 
\end{center}

\bigskip

We consider inverse periodic shadowing properties of discrete dynamical systems generated by diffeomorphisms of closed smooth manifolds.
We show that the $C^1$-interior of the set of all diffeomorphisms having so-called inverse periodic shadowing property coincides with the
set of $\Omega$-stable diffeomorphisms. The equivalence of Lipschitz inverse periodic shadowing property and hyperbolicity of
the closure of all periodic points is proved. Besides, we prove that the set of all diffeomorphisms that have Lipschitz inverse 
periodic shadowing property and whose periodic points are dense in the nonwandering set coincides with the set of Axiom A diffeomorphisms.

\bigskip
\bigskip

\section{Introduction. Main definitions. Statement of results.}

Theory of shadowing studies the problem of closeness of approximate and exact trajectories of dynamical systems on unbounded time intervals. Since the notions of an approximate trajectory and closeness can be formalized in several ways, various shadowing properties are considered (cf. \cite{1,2}). A problem about classical shadowing properties can be informally formulated in the following way: is it true that any sufficiently precise approximate trajectory of a dynamical system is close to some exact trajectory? In the paper \cite{4} an inverse problem was considered for the first time. This problem can be informally formulated as follows: suppose we have a method (pseudomethod) that generates approximate trajectories (pseudotrajectories) of a dynamical system; is it true that any exact trajectory is close to some pseudotrajectory generated by the pseudomethod?

In the present paper we consider discrete dynamical systems generated by $C^1$-diffeomorphisms of closed smooth manifolds. One of the main tasks of theory of shadowing is to characterise sets of all diffeomorphisms having some shadowing property or the $C^1$-interiors of such sets. One has to consider the $C^1$-interiors of the sets because for most cases it is difficult or maybe impossible to characterise this sets in terms of hyperbolic theory of dynamical systems. For many shadowing properties their $C^1$-interiors coincide with the set of structurally stable or $\Omega$-stable diffeomorphisms. K. Sakai (cf. \cite{5}) proved that the $C^1$-interior of the set of all diffeomorphisms having the standard shadowing property (also known as pseudo-orbit tracing property) coincides with the set of structurally stable diffeomorphisms. 

However, recently S.Yu.Pilyugin and S.B.Tikhomirov showed \cite{3} that Lipschitz shadowing property is equivalent to structural stability; and S.Yu. Pilyugin, S.B. Tikhomirov, and the author proved \cite{6} that so-called Lipschitz periodic shadowing property is equivalent to $\Omega$-stability, and the $C^1$-interior of the set of all diffeomorphisms having periodic shadowing property coincides with the set of $\Omega$-stable diffeomorphisms. Besides, recently S.Yu.Pilyugin, D.I. Todorov, and G.I.Wolfson (cf. \cite{7}) developed and improved technique from the paper \cite{3} to prove the equivalence of Lipschitz inverse shadowing property (for certain classes of pseudomethods) and structural stability.

While studying how well exact trajectories are approximated by pseudotrajectories generated by pseudomethods, it is not obligatory to consider all points of the phase space. It is possible to restrict research to some important invariant subset of the phase space, e.g., to the set of all periodic points. The corresponding shadowing properties can be called inverse periodic shadowing properties (rigorous definitions are given below). The present paper is devoted to study of such properties. First we give main definitions, and then formulate our results.

Let $M$ be a closed smooth manifold with Riemannian metric dist, $f$ be a diffeomorphism of the manifold $M$.

We say that a sequence of continuous mappings $\{\Psi_k\}_{k\in\mathbb{Z}}$ is a $d$-pseudomethod if
\begin{equation}
\label{pseudo}
\mbox{dist}(\Psi_k(x),f(x))\leq d\quad\textrm{for all } x\in M.
\end{equation}

We say that a sequence $\theta=\{x_k\}\subset M$ is a pseudotrajectory generated by a $d$-pseudomethod $\Psi=\{\Psi_k\}$ if
\begin{equation}
\label{pseudo2}
x_{k+1}=\Psi_k(x_k)\quad\forall k\in\mathbb{Z}.
\end{equation}

We say that a diffeomorphism $f$ has inverse periodic shadowing property InvPerSh if for any positive number $\epsilon$ there exists a positive number $d$ such that for any periodic point $p$ and for any $d$-pseudomethod $\Psi=\{\Psi_k\}$ there exists a pseudotrajectory $\theta=\{x_k\}$ generated by the $d$-pseudomethod $\Psi$ such that
\begin{equation}
\label{basis}
\mbox{dist}(x_k,f^k(p))<\epsilon\quad\forall k\in\mathbb{Z}.
\end{equation}
If inequalities $(\ref{basis})$ hold, we say that the pseudomethod $\{\Psi_k\}$ $\epsilon$-shadows the point $p$, or that the point $p$ is $\epsilon$-shadowed by the pseudomethod $\{\Psi_k\}$.

Let us define a Lipschitz version of InvPerSh. We say that a diffeomorphism $f$ has Lipschitz inverse periodic shadowing property LipInvPerSh if there exist positive numbers $L$ and $d_0$ such that for any periodic point $p$ and for any $d$-pseudomethod $\Psi=\{\Psi_k\}$ with $d\leq d_0$ there exists 
a pseudotrajectory $\theta=\{x_k\}$ generated by the $d$-pseudomethod $\Psi$ such that 
$$\mbox{dist}(x_k,f^k(p))\leq Ld\quad\forall k\in\mathbb{Z}.$$

As usual, we denote  by $\Omega S$ the set of all $\Omega$-stable diffeomorphisms. We denote by $\mbox{Int}^1(A)$ the $C^1$-interior of a set $A$ of diffeomorphisms. Finally we denote by $\mbox{Per}(f)$ the set of all periodic points of a diffeomorphism $f$. 

Our main result is the following theorem:

\textbf{Theorem.}
\begin{enumerate}
	\item[1)] $\mbox{Int}^1(\mbox{InvPerSh})=\Omega S$;
	
	\item[2)] LipInvPerSh is equivalent to hyperbolicity of the set $\mbox{Cl}(\mbox{Per}(f))$;

 \item[3)] if we denote  by LIPS the set of all diffeomorphisms that have LipInvPerSh and whose periodic points are dense in the nonwandering set, then LIPS coincides with the set of Axiom A diffeomorphisms.
\end{enumerate}

\textbf{Remark 1.} There are several ways of inroducing pseudomethods. One of them was described above. Such pseudomethods are called pseudomethods of class $\Theta_s$. The mappings $\Psi_k$ from the definition of pseudomethods of class $\Theta_s$ are close to the diffeomorphism $f$. It is possible to consider pseudomethods defined in the following way:

We say that a sequence $\{\Psi_k\}$ of continuous mappings is a $d$-pseudomethod (of class $\Theta_t$) if (instead of inequalities $(\ref{pseudo})$) the following inequalities hold:
$$\mbox{dist}(\Psi_{k+1}(x),f(\Psi_k(x)))\leq d\quad\textrm{for all }x\in M,\ k\in\mathbb{Z},$$
and we say that a sequence $\{x_k\}$ is a pseudotrajectory generated by the pseudomethod $\{\Psi_k\}$ if
$$x_k=\Psi_k(x_0)\quad\textrm{for all }k\in\mathbb{Z}.$$

Such pseudomethods are called pseudomethods of class $\Theta_t$. There are examples of pseudomethods of class $\Theta_s$ that do not belong to class $\Theta_t$ and vice versa. It is possible to introduce inverse periodic shadowing properties for pseudomethods of class $\Theta_t$ and to prove for them an analog of Theorem. The proof will be similar; however, the pseudomethods should be constructed in a different way. In particular, for any fixed $k$ the mappings $\Psi_k$ should be constant. We do not describe here the process of construction of such constant mappings, as they can be easily constructed basing on the pseudomethods of class $\Theta_s$ whose construction will be described in details.

\textbf{Scheme of the Theorem proof.} In essence, we use the strategy from the paper \cite{6}.

1) Denote by HP the set of all diffeomorphisms that do not have nonhyperbolic periodic points. Actually, it is proved in the paper \cite{9} that HP is a subset of LipInvPerSh. We use the result of Aoki and Hayashi (cf. \cite{10,11}) that states the equality $\mbox{Int}^1(\mbox{HP})=\Omega S$. Thus, in order to prove item 1) of Theorem, it is enough to establish the inclusion $\mbox{Int}^1(\mbox{InvPerSh})\subset\mbox{HP}$. In order to get this inclusion, we $C^1$-slightly perturb a diffeomorphism with a nonhyperbolic periodic point so that it does not have InvPerSh.

2) Let us describe the scheme of the proof of item 2). First we prove that LipInvPerSh implies hyperbolicity of any periodic point, then we prove uniform hyperbolicity of the set of all periodic points, and finally we establish the fact that the closure of all periodic points is a hyperbolic set.

3) In essence, item 2) implies item 3) because, by definition, Axiom A is equivalent to hyperbolicity of the nonwandering set and density of periodic points in the nonwandering set.

\section{Technical remarks and the exponential map}

The proof of Theorem consists of several lemmas. The proof of most of them uses a known technique based on the exponential map, which allows to transfer results from points in the manifold to tangent vectors and vice versa. We shall describe this technique in this section. Besides, it is convenient to introduce in this section main notations, which will be used in the sequel.  

Let $M$ be a closed smooth manifold with Riemannian metric dist. Denote by $\exp:TM\mapsto M$ the standard exponential map and by $\exp_x$ its restriction to $T_xM$, the tangent space at point $x$.

Denote by $N(r,x)$ the $r$-neighborhood of the point $x$ in the manifold $M$, by $B_T(r,y)$ the ball of radius $r$ with center at the point $y$ in the space $T_xM$. Besides, we denote by $B(r,A)$ the $r$-neighborhood of the set $A$ that is a subset of an Euclidean space.

There exists a positive number $r<1$ such that for any point $x\in M$ the mapping $\exp_x$ is a diffeomorphism of the set $B_T(r,0)$ onto its image, and the mapping $\exp^{-1}_x$ is a diffeomorphism of the set $N(r,x)$ onto its image. Besides, we assume that the number $r$ is chosen so small that the following holds:
\begin{equation}
\label{prop1}
\frac{\mbox{dist}(\exp_x(v),\exp_x(w))}{|v-w|}\leq 2\quad\mbox{for }v,w\in B_T(r,x),v\neq w;
\end{equation}
\begin{equation}
\label{prop2}
\frac{|\exp^{-1}_x(y)-\exp^{-1}_x(z)|}{\mbox{dist}(y,z)}\leq 2\quad\mbox{for }y,z\in N(r,x), y\neq z.
\end{equation}
We can always get this conditions, because
\begin{equation}
\label{prop}
D\exp_x(0) = \mbox{id}.
\end{equation}

Let $p$ be a point of a diffeomorphism $f$ of the manifold $M$, let $p_k~=f^k(p)$ and $A_k=Df(p_k)$ for all $k\in\mathbb{Z}$ (this notations will be used in the sequel). Consider the mappings
\begin{equation}
\label{3.2.1}
F_k=\exp^{-1}_{p_{k+1}}\circ f\circ\exp_{p_k}:T_{p_k}M\mapsto T_{p_{k+1}}M.
\end{equation}
By the standard property $(\ref{prop})$ of the exponential map, $DF_k(0)=A_k$. We can always represent $F_k(v)$ in the following form:
$$F_k(v)=A_kv+\phi_k(v),\quad\mbox{where }\frac{|\phi_k(v)|}{|v|}\longrightarrow0\ \mbox{as }|v|\rightarrow0.$$

We denote by $O(p,f)$ the trajectory of the point $p$ of the diffeomorphism~$f$, i.e., $O(p,f)=\{p_k=f^k(p)\mid k\in\mathbb{Z}\}$.

Besides, we need the following auxiliary statement, which will be used for construction of pseudomethods:

\textbf{Proposition.} Let $f:B(b,O(p,f))\mapsto\mathbb{R}^n$ be a $C^1$-smooth map, $p$ be a periodic point of the map $f$ of the fundamental period $m$, let sets $B(b,p_1),\ldots, B(b,p_m)$ be disjoint. 

1) Let $\epsilon<b/2$ be a small number. Assume that $f(x) = p_{k+1} + A_k(x-p_k)$ for $x\in B(b,p_k)$ ($1\leq k\leq m$), i.e., $f$ is the linear map. Let $d<\epsilon/2$ be an arbitrary sufficiently small number. Assume that we have constructed a continuous map $\psi: B(b,O(p,f))\mapsto\mathbb{R}^n$ such that
\begin{equation}
\label{cond1}
|\psi(x) - A_k(x-p_k) - p_{k+1}|\leq d\quad \textrm{for all }x\in B(b,p_k),\ 1\leq k\leq m.
\end{equation}
Then there exists a map $\Psi:B(b,O(p,f))\mapsto\mathbb{R}^n$ such that $|\Psi(x) - f(x)|\leq\allowbreak\leq d$ for all $x\in B(b,p_k)$, $1\leq k\leq m$,
the mappings $\psi$ and $\Psi$ coincide on the set $B(\epsilon/2,O(p,f))$, and the maps $\Psi$ and $f$ coincide on the set $B(b,O(p,f))\backslash B(\epsilon,O(p,f))$.

2) Let $C$ be an arbitrary large number, $d$ be an arbitrary sufficiently small number such that $Cd < b/2$. Let 
$$f(x) = p_{k+1} + A_k(x-p_k) + \phi_k(x)\quad\textrm{for }x\in B(b,p_k), 1\leq k\leq m,$$
where $|\phi_k(x)|\longrightarrow 0$ as $|x|\rightarrow0$. Assume that $d$ is so small that 
$$|\phi_k(x)|\leq d/2\quad\textrm{for }|x|\leq Cd.$$

Finally assume that we have constructed a continuous map $\psi: B(b,O(p,f))\mapsto\mathbb{R}^n$ such that condition $(\ref{cond1})$ holds with $d/2$ instead of $d$. Then there exists a map $\Psi:B(b,O(p,f))\mapsto\mathbb{R}^n$ such that $|\Psi(x) - f(x)|\leq d$ for all $x\in B(b,O(p,f))$,
the maps $\psi$ and $\Psi$ coincide on the set $B(Cd/2,O(p,f))$, and the maps $\Psi$ and $f$ coincide on the set $B(b,O(p,f))\backslash B(Cd,O(p,f))$.

\textbf{Proof.} Let us start from item 1). Choose a smooth monotonous function $\beta:[0,+\infty)\mapsto[0,1]$ such that $\beta(x) = 0$ for
$x\leq\epsilon/2$, $\beta(x) = 1$ for 
$x\geq\epsilon$.

Define the map $\Psi$ by the following formula:
$$\Psi(x) = (1-\beta(|x-p_k|))\psi(x) + \beta(|x-p_k|)(p_{k+1} + A_k(x - p_k))$$
for $x\in B(b,p_k)$, $1\leq k\leq m$. 
For $x\in B(b,p_k)$ we have the formula 
$$|\Psi(x) - f(x)| = |\Psi(x) - (\beta(|x-p_k|) + 1-\beta(|x-p_k|))f(x)|\leq (1-\beta(x))d + 0 \leq d.$$
Clearly, the map $\Psi$ is the desired.

Item 2) can be proved in the similar way.

\textbf{Remark 2.} The restricitions on the number $d$ do not depend on the map~$\psi$ if the condition $(\ref{cond1})$ holds.

\section{Proof of the main result.}

In the paper \cite{9} it is proved that if $A$ is a hyperbolic set, then there exist positive numbers $L$ and $d_0$ such that for any point $p\in A$, any number $d\leq d_0$, and any $d$-pseudomethod $\Psi=\{\Psi_k\}$ there exists a pseudotrajectory $\{x_k\}$ generated by the pseudomethod $\Psi$ such that the analog of relation $(\ref{basis})$ with $\epsilon=Ld$ holds. Thus, the following lemma holds:

\textbf{Lemma 1.} If the set $\mbox{Cl}(\mbox{Per}(f))$ is hyperbolic, then $f$ has LipInvPerSh.

\textbf{Corollary.} $\Omega S\subset\mbox{Int}^1(\mbox{InvPerSh})$.

\textbf{Lemma 2.} $\mbox{Int}^1(\mbox{InvPerSh})\subset\Omega S$.

\textbf{Proof.} By the lemma of Hayashi and Aoki (\cite{10,11}), $\mbox{Int}^1(\mbox{HP})=\Omega S$.
That is why it is enough to prove that $\mbox{Int}^1(\mbox{InvPerSh})\subset\mbox{HP}$.
Assume that the diffeomorphism $f\in \mbox{Int}^1(\mbox{InvPerSh})$ and does not belong to HP.
Thus, there exists a neighborhood $W$ of the diffeomorphism $f$ in the $C^1$-topology such that
$W\subset\mbox{InvPerSh}$ and $W\cap\mbox{HP}=\emptyset$. The diffeomorphism $f$ has a nonhyperbolic periodic point $p$ of the fundamental period $m$, i.e., the operator $Df^m(p)$ has an eigenvalue $|\lambda|=1$. Without loss of generality, we assume that $\lambda$ is a purely complex number. The case of a real $\lambda$ can be treated in the similar way. 

At first we $C^1$-slightly perturb the diffeomorphism $f$ to get a diffeomorphism $h$ with certain properties that is linear in a neighborhood of its periodic trajectory $p_1,\ldots,p_m$. 

There exists a number $a\in (0,r)$ (recall that the number $r$ was defined in section 2) and a diffeomorphism $h\in W$ such that $h(p_j)=p_{j+1}$;
the point $p_j$ is assigned to 0 in coordinates $v_j=(\rho_j\cos\theta_j,\rho_j\sin\theta_j,w_j)_j$ in the space $T_{p_j}M$; and~if 
$$H_j=\exp^{-1}_{p_{j+1}}\circ h\circ\exp_{p_j}$$
 and $|v_j|\leq a$, then for some real number $\chi$ and natural number $\nu$ such that $\cos\nu\chi=1$ the following holds:
 $$H_j(v_j)=A_jv_j=(r_j\rho_j\cos(\theta_j+\chi),r_j\rho_j\sin(\theta_j+\chi),B_jw_j)_{j+1}$$
(where $B_j$ is a matrix of size $(n-2)\times(n-2)$, $n=\dim M$, $B_{j+m}=B_j$)
and $$r_0r_1\cdots r_{m-1}=1,\quad r_{j+m}=r_j.$$
Hereinafter we use the index $j$ after the brackets to emphasize that the vector is represented in the coordinates in the tangent space at the point $p_j$. 

Thus, the diffeomorphism $h$ is $C^1$-close to the diffeomorphism $f$, and the operator $Dh^m(p_0)$ has an eigenvalue $\lambda$ that is a root of degree $\nu$ of 1 and corresponds to a Jordan block of dimension one.

Choose a number $\bar{a}<a$ such that $B_T(\bar{a},0)_j\subset H_j^{-1}(B_T(a,0)_{j+1})$ for all $j$. We use the index $j$ after the brackets to emphasize that we work with a ball in the space $T_{p_j}M$.

Put 
$$R=2\max(r_0,\ldots,r_{m-1}).$$
We assume that the numbers $a$ and $\bar{a}$ were chosen so small that the neighborhoods $\exp_{p_k}(B_T(\bar{a},0)_k)$ are disjoint for all $1\leq k\leq m$.
Let $\epsilon_0 = \bar{a}/3$, $\epsilon=\epsilon_0/10$, and let $m\nu d<\epsilon/3$ be an arbitrary sufficiently small number.  

Define maps $\psi_k: \bigcup_{1\leq l\leq m}B_T(\bar{a},0)_l\mapsto \bigcup_{1\leq l\leq m}B_T(a,0)_l$ in the following way:
$$\psi_{k+m\nu} = \psi_k\quad\textrm{for all }k\in\mathbb{Z};$$
$$\psi_k(y)=A_ky+(dr_0\cdots r_k(\cos (k+1)\chi)/(2R^m),dr_0\cdots r_k(\sin (k+1)\chi)/(2R^m),0)_{k+1}$$
 for $0\leq k\leq m\nu - 1$, $y\in B_T(\bar{a},0)_k$; and $\psi_k(y) = H_l(y)$ for $y\in B_T(\bar{a},0)_l$, $l\neq k$.

The maps $\psi_k$ can be considered as mappings from the disjoint union of $n$-dimensional balls in $\mathbb{R}^n$ to the the disjoint union of larger $n$-dimensional balls in $\mathbb{R}^n$. Choose an arbitrary number $k\in\mathbb{Z}$. We observe that the maps $\psi_k$ satisfy condition $(\ref{cond1})$.
Thus, in essence, all conditions of Proposition from section 2 are satisfied. Since we can decrease $d$, we can assume that $d$ is a number from item 1) of Proposition that is applied to the number $\epsilon_0$ (not number $\epsilon$) as $\epsilon$ and the map $\psi_k$ (we put $b=\bar{a}$). By Remark 2, the choice of $d$ does not depend on the map $\psi_k$ if condition $(\ref{cond1})$ holds. Consequently, the map $\psi_k$ can be extended to the map $\Phi_k$, which coincides with $H_k$ on the sets $\bigcup_{1\leq l\leq m} B_T(\bar{a},p_l)  \backslash\bigcup_{1\leq l\leq m} B_T(\epsilon_0/2,p_l)$.
Put $\Psi_k(x) = \exp_{p_{l+1}}\circ \Phi_k\circ\exp^{-1}_{p_l}(x)$ for $x\in \exp_{p_l}(B_T(\bar{a},0)_l)$. If we set $\Psi_k = h$ on the complement of $\bigcup_{1\leq l\leq m}\exp_{p_l}B_T(\bar{a},0)_{l}$ in the manifold $M$, then $\Psi_k$ will be a continuous map defined on $M$. By property $(\ref{prop1})$, the maps $\{\Psi_k\}$ are a $2d$-pseudomethod.

However, no pseudotrajectory generated by this $2d$-pseudomethod can $\epsilon$-shadow the point $p$. In fact, suppose the contrary, consider a pseudotrajectory $\{y_k\}$ of the point $y$ defined by the analog of equalities $(\ref{pseudo2})$, and assume that this pseudotrajectory $\epsilon$-shadows the point $p$. Put $q=\exp_p^{-1}(y)=(q_1,q_2,q_3)$ (in the proof of Lemma 2 the first two components in such representation have dimension 1). Put $q^k = \exp^{-1}_{p_{k}}y_k $ for all $k\in\mathbb{Z}$, then $q^{k+1} = \Phi_kq^k$. Note that, by property $(\ref{prop2})$, the points $q^k$ are $2\epsilon$-close to 0 in the coordinates in $T_{p_k}M$. Let us emphasize that $2\epsilon<\epsilon_0/2$, and the maps $\Phi_k$ coincide with the maps $\psi_k$ in the $\epsilon_0/2$-neighborhood of zero (in the coordinates in $T_{p_k}M$). Hence,
$$q^{m\nu k}=A^{m\nu k}q+(km\nu d/(2R^m),0,0).$$ 
We denote by $pr_{1,2}$ the projection of a vector onto its first and second components. Then
\begin{equation}
\label{mage}
|pr_{1,2}q^{mk\nu}|\geq mk\nu d/(2R^m)-|pr_{1,2}q|^{mk\nu}=mk\nu d/(2R^m) - |pr_{1,2}q|.
\end{equation} 
By our assumptions, $|q_k|\leq2\epsilon$ for all $k\in\mathbb{Z}$.
However, by estimates $(\ref{mage})$, numbers $|\mbox{pr}_{1,2}q_k|$ will be larger than $3\epsilon$ for sufficiently large $k$. Hence, our assumptions are wrong, the diffeomorphism $f$ does not have any nonhyperbolic periodic points, and the $C^1$-interior of InvPerSh is contained in $\Omega S$.   

\textbf{Lemma 3.} If a diffeomorphism $f$ has LipInvPerSh, then any periodic point is hyperbolic.

\textbf{Proof.} Without loss of generality, we assume that the number $L$ from the definition of LipInvPerSh is natural.

Let $p$ be a nonhyperbolic point of the diffeomorphism $f$. We assume that $p$ is a fixed point, in order to simplify the notations. The general case can be treated in the similar way. In the case of a fixed point, the map $(\ref{3.2.1})$ is represented in the following way:
$$F(v)=\exp_p^{-1}\circ f\circ\exp_p(v)=Av+\phi(v),$$
and the matrix $A$ has an eigenvalue $|\lambda|=1$. We assume that $\lambda$ is a purely complex number. The case of a real $\lambda$ can be treated in the similar way. 

By the choice of coordinates, we assume that the matrix $A$ is represented in the form $\mbox{diag}(H_1,H_2)$, where
$$H_1=\left(\begin{matrix}Q&I&\cdots&\cdots& 0\\
0&Q&I&\cdots& 0\\
\cdots&\cdots&\cdots&\cdots&\cdots\\
0&0& \cdots&\cdots&Q
\end{matrix}\right)$$ 
and $$Q=\left(\begin{matrix}\cos\theta&\sin\theta\\-\sin\theta&\cos\theta\end{matrix}\right).$$
When the new coordinates were chosen, the numbers $L,d_0,r$ could change (recall that the number $r$ was defined in section 2). We denote by the same symbols the new constants for convenience. If $v$ is a two-dimensional vector, then $$|Qv|=|v|.$$
Let $2l$ be the dimension of the Jordan block $H_1$, and put $n=\dim M$.

Let us introduce some notations. Let $v$ be an $n$-dimensional vector. We denote by $pr_{i,j}v$ the two-dimensional vector that consists of the $i$th and $j$th components of the vector $v$. Let $V$ be a matrix of size $2\times 2$, and $W$ be a matrix of size $n\times n$. When we write $W = (0,V,0)^{k-1,k}$, we assume that the matrix $V$ is contained in the matrix $W$ on positions $(k-1,k)\times (k-1,k)$, and all other elements of the matrix $W$ are equal to zero.  

Choose a number $\bar{r}<r/2$ such that $B_T(\bar{r},0)\subset F^{-1}(B_T(r,0))$. Let $d$ be a sufficiently small arbitrary number such that $20Ld<\bar{r}/10$. 

Define the maps $\psi_k:B_T(\bar{r},0)\mapsto B_T(r,0)$ in the following way:
$$\psi_k(y) = Ay + (d/2)(0,Q^kw,0)^{(2l-1,2l)},$$
where $w$ is a two-dimensional unit vector.

Fix an arbitrary number $k$. The maps $\psi_k$ satisfy condition $(\ref{cond1})$.
We assume that the number $d$ is less than the corresponding number from item 2) of Proposition applied to $C=20L$, $b=\bar{r}$.
Let us emphasize that, by Remark 2, the choice of $d$ does not depend on the index $k$, since the maps $\psi_k$ satisfy the analog of condition $(\ref{cond1})$. Let $\Phi_k$ be the analog of the map $\Psi_k$ constructed in item 2) of Proposition. The map $\Phi_k$ coincides with $F$ on the set $B_T(\bar{r},0)\backslash B_T(10Ld,0)$. Hence, the map $\Psi_k = \exp_p\circ\Phi_k\circ\exp^{-1}_p$  can be extended to the continuous map $\Psi_k$ of the manifold $M$. By property $(\ref{prop1})$, the maps $\{\Psi_k\}$ are a $2d$-pseudomethod.

By construction, for $k>l$
$$
|\Phi_{k-1}\circ\ldots\circ\Phi_0(q)| =  |A^kq + c_k^1d(0,Q^{k-l+1}w,0)^{(1,2)} + 
\ldots +
$$
$$
+ c_k^{l-1}d(0,Q^{k-1}w,0)^{(2l-3,2l-2)}
 +
d(0,kQ^kw,0)^{(2l-1,2l)})|,
$$
where $c_k^m$ are some positive numbers. The following inequality holds:
\begin{equation}
\label{form}
|pr_{(2l-1,2l)}(\Phi_{k-1}\circ\ldots\circ\Phi_0(q))|\geq
kd - |pr_{(2l-1,2l)}(A^kq)|\geq kd - |pr_{(2l-1,2l)}q|^k.
\end{equation}

Assume that the point $p$ is $2Ld$-shadowed by a pseudotrajectory $\{y_k\}$ of some point $y$ generated by the pseudomethod $\{\Psi_k\}$. Put $q_k = \exp^{-1}_{p_{k}}y_k $ for all $k\in\mathbb{Z}$, then $q_{k+1} = \Phi_kq_k$.
We observe that, by inequalities~$(\ref{prop2})$, 
\begin{equation}
\label{contr1}
|pr_{(2l-1,2l)}q_k|\leq|q_k|\leq 4Ld\quad\textrm{for all }k\in\mathbb{Z}
\end{equation} 
(in the coordinates in $T_{p_k}M$). Note that the maps $\Phi_k$ and $\psi_k$ coincide on the $10Ld$-neighborhood of zero.

By relations $(\ref{form})$ and $(\ref{contr1})$,
$$|pr_{(2l-1,2l)}q_{10L}|\geq (10Ld - |pr_{(2l-1,2l)}q|)\geq 6Ld.$$
The last inequality contradicts to inequality $(\ref{contr1})$. Hence, the diffeomorphism $f$ does not have any nonhyperbolic periodic points.

\textbf{Lemma 4.} If a diffeomorphism $f$ has LipInvPerSh, then all periodic points of the diffeomorphism $f$ are uniformly hyperbolic (i.e. they are hyperbolic with the same constants $C$ and $\lambda$). In other words, there exist constants $C>0$ and $0<\lambda<1$ depending only on $L$ such that for any periodic point $p$ of the diffeomorphism $f$ there exist $Df$-invariant complementary subspaces $S(p)$ and $U(p)$ of the tangent space $T_{p_k}M$ such that 
$$|Df^j(p)v|\leq C\lambda^j|v|\quad\textrm{for }v\in S(p),\ j\geq0;$$
$$|Df^{-j}(p)v|\leq C\lambda^j|v|\quad\textrm{for }v\in U(p),\ j\geq0.$$

\textbf{Proof.} Without loss of generality, we assume that the number $L$ from the property  LipInvPerSh is natural.

Let $p$ be a periodic point of a period $m$. Denote by $m_0$ the fundamental period of the point $p$. Put $p_i=f^i(p)$, $A_i=Df(p_i)$, and $B=Df^m(p)$.
By Lemma~3, the point $p$ is a hyperbolic periodic point. Hence,  there exist complementary $Df$-invariant linear spaces $S(p)$ and $U(p)$ at point $p$, and this spaces satisfy the conditions
\begin{equation}
\label{3.4.13}
\lim_{n\rightarrow+\infty}B^nv_s=\lim_{n\rightarrow+\infty}B^{-n}v_u=0\quad\textrm{for } v_s\in S(p), v_u\in U(p).
\end{equation} 
Consider an arbitrary nonzero vector $v_u\in U(p)$. Put $e_0=v_u/|v_u|$. Consider the sequence
$$a_0=\tau,\quad a_{i+1} = a_i|A_ie_i| - 1,$$
where $e_{i+1} = A_ie_i/|A_ie_i|$,
and the number $\tau$ is chosen such that $a_m=0$. An explicit formula for the number $\tau$ that satisfies all the required conditions is given in the paper \cite{6}. Note that $a_j>0$ for all $0\leq j\leq m-1$ (since the inequality $a_j\leq 0$ implies the inequality $a_{j+1}<0$, whereas $a_m=0$).

It follows from relations $(\ref{3.4.13})$ that there exists a number $n>0$ such that
\begin{equation}
\label{3.4.16}
|B^{-n}\tau e_0|<1.
\end{equation}

Consider a finite sequence $w_i\in T_{p_i}M$ for $0\leq i\leq m(n+1)$ given by the equalities
$$w_i=a_ie_i\quad\textrm{for } i\in\{0,\ldots,m-1\},$$
$$w_m=B^{-n}\tau e_0,$$
$$w_{m+1+i}=A_iw_{m+i}\quad\textrm{for } i\in\{0,\ldots,mn-1\}.$$
Note that
$$w_{km}=B^{k-1-n}\tau e_0\quad\textrm{for } k\in\{1,\ldots,n+1\}.$$
Thus, we can consider the sequence $\{w_i\}$ as an $m(n+1)$-periodic sequence that is well-defined for all $i\in\mathbb{Z}$.

Put $N>\max_{k\in\mathbb{Z}}|w_k|$. Since we can increase the number $N$ if necessary, we assume that $N>20L$. It is clear that if all vectors from the sequence $\{w_k\}$ are multiplied by $d$, then the maximum will increase in $d$ times too. 

Choose a number $\epsilon_1<r$ such that the neighborhoods $\exp_{p_j}B_T(\epsilon_1,0)_j$ are disjoint for all $1\leq j\leq m_0$. Let us emphasize that the index $j$ after the brackets means that we work with a ball in the space $T_{p_j}M$. Choose a number $\epsilon_2\leq\epsilon_1$ such that $B_T(\epsilon_2,0)_j\subset F_j^{-1}(B_T(\epsilon_1,0))_j$ for all $1\leq j\leq m_0$.  

We assume that the number $d$ is a sufficiently small arbitrary number such that $100Nd<\epsilon_2$.

Define the maps $\psi_k:\bigcup_{1\leq i\leq m_0} B_T(\epsilon_2,0)_i\mapsto\bigcup_{1\leq i\leq m_0} B_T(\epsilon_1,0)_i$ in the following way:

\begin{enumerate}
\item[1)] $\psi_{k + m(n+1)} = \psi_k$ for all $k\in\mathbb{Z}$;

\item[2)] if $y\in B_T(\epsilon_2,0)_k$, then
\begin{enumerate}
\item[2.1)] $\psi_k(y)=A_ky - de_{k+1}$ for $0\leq k \leq m-2$,
\item[2.2)] $\psi_{m-1}(y)=A_{m-1}y - de_{m} + B^{-n}\tau de_{0}$ for $k=m-1$,
\item[2.3)] $\psi_k(y)=A_ky$ for $m\leq k\leq mn+m-1$;
\end{enumerate}

\item[3)] for other $y$ ($y\in B_T(\epsilon_2,0)_l$ and $l-k$ is not a multiple of $m_0$, the fundamental period of the point $p$) $\psi_k(y)=F_l(y)$.
\end{enumerate}

Let us show that the following equalities hold:
\begin{equation}
\label{defW1}
\psi_{k-1}\circ\ldots\circ\psi_0(w_0d)=w_kd\quad\textrm{for all }k\geq 1,
\end{equation}
\begin{equation}
\label{defW2}
\psi^{-1}_k\circ\ldots\circ\psi^{-1}_{-1}(w_0d)=w_kd\quad\textrm{for all }k\leq 0.
\end{equation}

Indeed, for any $0\leq k\leq m-2$
$$dw_{k+1} = da_{k+1}e_{k+1} =  d((a_k|A_ke_k| - 1)/|A_ke_k|)A_ke_k = A_kdw_k - de_{k+1} = \psi_k(dw_k).$$
Since, by the choice of $\tau$, $1=a_{m-1}|A_{m-1}e_{m-1}| = |w_{m-1}||A_{m-1}e_{m-1}|=\allowbreak=|A_{m-1}w_{m-1}|$, we have $A_{m-1}w_{m-1} = e_{m}$, hence,
$$\psi_{m-1}(dw_{m-1}) = A_{m-1}dw_{m-1} - de_{m} + B^{-n}\tau de_{0} = B^{-n}\tau de_{0} = dw_m.$$
Thus, we proved equalities $(\ref{defW1})$ and $(\ref{defW2})$. Note that the maps $\psi_k$ were constructed in such way that this equalities would hold.

Observe that, by inequality $(\ref{3.4.16})$, the maps $\psi_k$ satisfy the analog of condition $(\ref{cond1})$ with $2d$ instead of $d/2$. Fix an arbitrary $k\in\mathbb{Z}$. We see that all conditions of item 2) of Proposition are satisfied. We apply Proposition to $C = 100N$, $b=\epsilon_2$. Since we can decrease the number $d$, we can assume that it is less than $d/4$, where $d$ is from item 2) of Proposition. Let us emphasize that, by Remark 2, the number $d$ does not depend on the index~$k$. Denote by $\Phi_k$ the analog of the map $\Psi_k$ constructed in Proposition. By the statement of Proposition, $|\Phi_k(x) - F_j(x)|\leq 4d$ for $x\in\allowbreak B_T(\epsilon_2,0)_j$ and $1\leq j\leq m_0$. Let us emphasize that the map $\Phi_k$ coincides with the map $F_k$ on the set $\bigcup_{0\leq j\leq m_0} B_T(\epsilon_2,0)_j\backslash \bigcup_{0\leq j\leq m_0} B_T(50Nd,0)_j$. Consider the maps $\Psi_k$ given by the formula
$$\Psi_k(y) = \exp_{p_{l+1}}\circ F_k\circ\exp^{-1}_{p_l}(y)\quad\textrm{for }y\in \exp_{p_l}(B_T(\epsilon_2,0)_l),\ 1\leq l\leq m_0.$$
Clearly, if we define the maps $\Psi_k$ to be equal to the diffeomorphism $f$ for all other points of the manifold $M$, then the maps $\Psi_k$ will remain being continuous maps. Note that, by inequalities $(\ref{prop1})$, the maps $\Psi_k$ generate an $8d$-pseudomethod.

By our assumptions, the point $p$ is $8Ld$-shadowed by one of pseudotrajectories generated by the pseudomethod $\Psi =\{\Psi_k\}$. It follows that there exists a pseudotrajectory $\{y_k\}$ of the point $y$ defined by the analog of equalities $(\ref{pseudo2})$ such that the point $p$ is $8Ld$-shadowed by this pseudotrajectory. It follows from our assumptions that
\begin{equation}
\label{why}
|q_k|\leq16Ld.
\end{equation}
Note that $16Ld<50Ld$, and the maps $\Phi_k$ and $\psi_k$ coincide on the $50Ld$-neighborhood of zero (in the space $T_{p_k}M$).

It is easily seen that
$$\Phi_0(q-w_0d + w_0d)=A_0(q-w_0d) + \Phi_0(w_0d)=A_0(q-w_0d)+w_1d,$$
$$\Phi_{k-1}\circ\cdots\circ\Phi_0(q-w_0d + w_0d)=A_{k-1}\cdots A_0(q-w_0d) + \Phi_{k-1}\circ\cdots\circ\Phi_0(w_0d)=$$
\begin{equation}
\label{formula1}
=A_{k-1}\cdots A_0(q-w_0d)+w_kd,\quad k\geq1,
\end{equation}
\begin{equation}
\label{formula2}
\Phi^{-1}_{-k}\circ\cdots\circ\Phi^{-1}_{-1}(q-w_0d + w_0d)=A^{-1}_{-k}\cdots A^{-1}_{-1}(q-w_0d)+w_{-k}d,\quad k\geq1.
\end{equation}

In $(\ref{formula1})$ and $(\ref{formula2})$ the norm of the second term is estimated from above by $Nd>20Ld$, whereas, by hyperbolicity of the point $p$, the first term in one of this formulae has a large norm (larger than $2Nd$) for $k$ with large absolute values (i.e. the norm of the first term is much larger than the norm of the second term) if $q\neq w_0d$.
Thus, the point $p$ can be shadowed only by the pseudotrajectory corresponding to the vector $w_0d$, i.e. $q=w_0d$. But then inequalities $(\ref{why})$ imply the estimates
$$|a_k|=|w_k|\leq 16L\quad\mbox{for }0\leq k\leq m-1.$$
This estimates imply the desired hyperbolicity estimates for any preliminary fixed vector $v$ (cf. \cite{6} for detailed explanation).
\\

Note that, actually, Lemma 4 implies item 2) of Theorem (and, consequently, item 3) too); since it is proved in the paper \cite{6} that if the statement of Lemma 4 holds, i.e., the set $\mbox{Per}(f)$ has a hyperbolic structure, then the set $\mbox{Cl}(\mbox{Per}(f))$ has a hyperbolic structure too, i.e., the set $\mbox{Cl}(\mbox{Per}(f))$ is hyperbolic.

\end{document}